\documentclass[12pt]{amsart}

\usepackage{amssymb,amsfonts}
\usepackage[all,arc]{xy}
\usepackage{enumerate}
\usepackage{mathrsfs}
\usepackage{graphicx}

\usepackage{mathtools}

\usepackage{yfonts}

\newtheorem{theorem}{Theorem}[section]
\newtheorem{corollary}[theorem]{Corollary}
\newtheorem{proposition}[theorem]{Proposition}
\newtheorem{lemma}[theorem]{Lemma}

\theoremstyle{definition}
\newtheorem{definition}[theorem]{Definition}

\newtheorem{example}[theorem]{Example}

\newtheorem{question}[theorem]{Question}
\newtheorem{remark}[theorem]{Remark}

\DeclareMathOperator{\Max}{Max}
\DeclareMathOperator{\Spec}{Spec}
\DeclareMathOperator{\Id}{Id}

\numberwithin{equation}{section}

\begin{document}
	\title{Some remarks on the comparability of ideals in semirings}
	
	\author{H. Behzadipour}
	
	\address{School of Electrical and Computer Engineering\\ University College of Engineering\\ University of Tehran\\ Tehran\\ Iran}
	
	\email{hussein.behzadipour@gmail.com}
	
	\author{P. Nasehpour}
	
	\address{Department of Basic Sciences\\ Golpayegan College of Engineering\\ Isfahan University of
		Technology\\ Golpayegan 87717-67498\\ Iran.}
	\email{p.nasehpour@iut.ac.ir, nasehpour@gmail.com}
	
	
	
\begin{abstract}
A semiring is uniserial if its ideals are totally ordered by inclusion. First, we show that a semiring $S$ is uniserial if and only if the matrix semiring $M_n(S)$ is uniserial. As a generalization of valuation semirings, we also investigate those semirings whose prime ideals are linearly ordered by inclusion. For example, we prove that the prime ideals of a commutative semiring $S$ are linearly ordered if and only if for each $x,y \in S$, there is a positive integer $n$ such that either $x|y^n$ or $y|x^n$. Then, we introduce and characterize pseudo-valuation semidomains. It is shown that prime ideals of pseudo-valuation semidomains and also of the divided ones are linearly ordered.
\end{abstract}
	
	\maketitle
	
	\section{Introduction}
	
	Hemirings and semirings are ring-like algebraic structures in which subtraction is either impossible or disallowed. Since the language of semiring theory is not standardized yet, we begin our paper by explaining what we mean by a hemiring and a semiring.
	
	A hemiring is an algebraic structure $(H,+,\cdot,0)$ satisfying the following axioms:
	
	\begin{enumerate}
	\item $(H,+,0)$ is a commutative monoid,
	\item $(H,\cdot)$ is a semiring,
	\item $a(b+c) = ab+ac$ and $(b+c)a = ba + ca$, for all $a,b,c \in H$,
	\item $a\cdot 0 = 0 = 0 \cdot a$ for all $a\in H$. 
\end{enumerate}

    In \S\ref{sec:orderlrannihilators}, we investigate the finiteness or linearity conditions \cite{DaveyPriestley2002} of the poset of some special subsets of a hemiring. More precisely, if $X$ is a nonempty subset of a hemiring $H$, the left (right) annihilator of $X$ is defined as follows: \[\textswab{l}(X)= \{h\in H: hX = \{0\}\}~(\textswab{r}(X)= \{h\in H: Xh = \{0\}\}).\] In Corollary \ref{orderlrannihilators}, we prove that in an arbitrary hemiring $H$, the ascending chain condition for right annihilators is equivalent to the descending chain condition for left annihilators. We also prove that right annihilators are linearly ordered by inclusion if and only if left annihilators are so.
    
    An algebraic structure $(S,+,\cdot,0,1)$ is a semiring if $(S,+,\cdot,0)$ is a hemiring , $1\neq 0$, and $s1 = 1s = s$ for all $s\in S$. A nonempty subset $I$ of a semiring $S$ is said to be an ideal of $S$ if $x+y \in S$ and $sx , xs \in S$ for all $x,y \in I$ and $s\in S$ \cite{Bourne1951}. The ideal $I$ is proper if $I \neq S$. The family of ideals of a semiring $S$ is denoted by $\Id(S)$. 
    
    A proper ideal $P$ of a semiring $S$ is prime if $IJ \subseteq P$ implies either $I \subseteq P$ or $J \subseteq P$ for all ideals $I$ and $J$ of $S$. A semiring $S$ is commutative if $ab = ba$, for all $a,b \in S$. In a commutative semiring $S$, $P$ is a prime ideal of $S$ if and only if $ab\in P$ implies either $a\in P$ or $b \in P$, for all $a,b \in S$ \cite[Corollary 7.6]{Golan1999(b)}. The family of prime ideals of a semiring $S$ is denoted by $\Spec(S)$. For more on the ideals of commutative semirings, see \cite{Nasehpour2020,NasehpourQ2018}.
    
     Note that a ring $R$ is said to be uniserial if the family of its (two-sided) ideals is a chain \cite{Puninski2001}. Similarly, we define a semiring $S$ to be uniserial if its ideals are linearly ordered by inclusion. In \S\ref{sec:chainconditionsmatrixsemirings}, we investigate the chain conditions of the ideals of matrix semirings. In particular, we show that if $S$ is a semiring, then the following statements hold:
    
    \begin{enumerate}
    	\item The semiring $S$ is uniserial if and only if the matrix semiring $M_n(S)$ is uniserial.
    	\item The prime ideals of the semiring $S$ are linearly ordered by inclusion if and only if the prime ideals of $M_n(S)$ are so.
    	\item The subtractive (prime) ideals of $S$ are linearly ordered by inclusion if and only if the subtractive (prime) ideals of $M_n(S)$ are so.
    \end{enumerate}

    Let us recall that a proper ideal $M$ of a semiring $S$ is maximal if there are no other ideals to be properly between $M$ and $S$. The family of maximal ideals of a semiring $S$ is denoted by $\Max(S)$. 

    In \S\ref{sec:pairsofsemiringswiththesameprimeideals}, we introduce the concept of pseudo-valuation subsemiring (check Definition \ref{pseudovaluationsubsemiringdef}) which is a generalization of its counterpart in rings (see \cite{AyacheEchi2006} and \cite{KumarGaur2020}). By definition, a commutative semiring $S$ is a semidomain if it is multiplicatively cancellative, i.e. $ab =ac$ implies $b=c$ for all $a,b,c \in S$ with $a\neq 0$. Now, let $S \subset U$ be semidomains. We say $S$ is a pseudo-valuation subsemiring of $U$ if for each $u \in U \setminus S$ and non-unit $s\in S$, we have $u^{-1} s \in S$. In Theorem \ref{pseudo-valuationsubsemiringthm}, we prove that if $S \subseteq T \subset U$ are semidomains with $\Max(T) \subseteq \Spec(S)$ and $\Max(S) \subseteq \Spec(T)$. Then, $S$ is a pseudo-valuation subsemiring of $U$ if and only if $T$ is so. A corollary to this result is that if $S \subseteq T \subset U$ are semidomains with $\Spec(S) = \Spec(T)$, then $S$ is a pseudo-valuation subsemiring of $U$ if and only if $T$ is so. Note that this is a generalization of Theorem 2.2 in \cite{KumarGaur2020}.    
	
	In \S\ref{sec:lopi}, we continue our investigations of the chain conditions of commutative semirings and divided semidomains which we had discussed in \cite{BehzadipourNasehpour2020}. 
	
	Let us recall that a semidomain is a valuation semiring if and only if the poset of its ideals is a chain, i.e., its ideals are linearly ordered by inclusion (see Theorem 2.4 in \cite{Nasehpour2018}). Also note that if $a$ and $b$ are elements of a commutative semiring $S$, it is said that $b$ divides $a$ or is a divisor of $a$ and written $b \mid a$, if there exists an element $x \in S$ such that $a = bx$. It is evident that this is equivalent to say that $(a) \subseteq (b)$, where by $(a)$, we mean the principal ideal of $S$ generated by $a$. As a generalization of valuation semirings, we discuss those semirings that their prime ideals are linearly ordered by inclusion. Similar to \cite{Badawi1995}, in Theorem \ref{LOPmain}, we prove that the prime ideals of a commutative semiring $S$ are linearly ordered if and only if for each $x,y \in S$, there is a positive integer $n$ such that either $x|y^n$ or $y|x^n$. It is worth mentioning that with the help of this statement, in Theorem \ref{GCDvaluationvasconcelos}, we show that a semidomain is a valuation semiring if and only if it is a GCD semidomain which has prime ideals that are linearly ordered. This is exactly the semiring version of a result by Vasconcelos in \cite{Vasconcelos1972}, though our proof is similar to Badawi's proof in \cite{Badawi1995} which is different from the proof given by Vasconcelos.
	
	Using techniques adapted from ring theory, it is straightforward to show that if $S$ is a semidomain, then $S$ can be embedded in a semifield, called its semifield of fractions, denoted by $F(S)$ (see p. 22 in \cite{Golan1999(a)}). In fact, the localization of the $S$ at $S \setminus\{0\}$ gives the semifield $F(S)$ (for more details on this, see \S11 in \cite{Golan1999(b)}).  
	
	Based on the papers \cite{Badawi1995,HedstromHouston1978}, we introduce pseudo-valuation semidomains as follows (see Definition \ref{PVS}): We define a prime ideal $\mathfrak{p}$ of a commutative semiring $S$ to be strongly prime if whenever $x,y \in F(S)$ and $xy \in \mathfrak{p}$, then $x \in \mathfrak{p}$ or $y \in \mathfrak{p}$. We say that a semidomain $S$ is a pseudo-valuation semidomain (abbreviated as PVS) if every prime ideal of $S$ is strongly prime. Then we prove that pseudo-valuation semidomains are examples of those semirings that have linearly ordered prime ideals (check Proposition \ref{PVSLOP}). 
	
	In Theorem \ref{PVSchar1}, we show that a semidomain $S$ is a PVS if and only if $S$ is quasi-local and the only maximal ideal $\mathfrak{m}$ of the commutative semiring $S$ is strongly prime. In Theorem \ref{subsetlocal}, which is the semiring version of Proposition 3.1 in \cite{Anderson1979}, we show that if $(S,\mathfrak{m})$ is a quasi-local semidomain, then, $S$ is a PVS if and only if for every $x \in F(S)$, either $xS \subseteq \mathfrak{m}$ or $\mathfrak{m} \subseteq xS$. It is worth mentioning that in Theorem \ref{PVSchar2}, we give another characterization of pseudo-valuation semidomains.
	
	In Theorem \ref{LOP}, we show that if $(S,\mathfrak{m})$ is a quasi-local semidomain and for each $x \in F(S)$, $x \mathfrak{m} \subseteq \mathfrak{m}$ or $\mathfrak{m} \subseteq x \mathfrak{m}$, then the set of prime ideals of $S$ is linearly ordered.
	
	Another family of commutative semirings which have linearly ordered prime ideals is the family of divided semidomains introduced in \cite{BehzadipourNasehpour2020} (see Corollary \ref{dividedLO}). Finally, in Theorem \ref{dividedchar}, we prove that a semidomain $S$ is divided if and only if for each $x,y \in S$, either $x|y$ or $y|x^n$  for some $n \geq 1$.
	
	Ring-like algebraic structures such as hemirings and semirings have many applications in computer science and engineering \cite[p. 225]{Golan2003}. Also, note that semirings are considered to be interesting generalizations of bounded distributive lattices and rings \cite[Example 1.5]{Golan1999(b)}. For more on semirings and their applications, one may refer to \cite{Glazek2002,Golan1999(a),Golan1999(b),Golan2003,GondranMinoux2008,HebischWeinert1998,KuichSalomaa1986,NasehpourHistory}.
	
	\section{Chain conditions on left and right annihilators of hemirings}\label{sec:orderlrannihilators}
	
	Let $(P_1, \leq_1)$ and $(P_2, \leq_2)$ be posets. A function $f: P_1 \rightarrow P_2$ is called to be a poset antihomomorphism if it is oreder-reversing, i.e. $x \leq_1 y$ implies $f(y) \leq_2 f(x)$, for all $x,y \in P_1$. A poset antihomomorphism from a poset into itself is called to be a poset antiendomorphism. 
	
	\begin{definition}
    Let $(P,\leq)$ be a poset and $f$ a poset antiendomorphism. We say an element $y\in P$ is an $f$-element of $P$ if there is an element $x$ in $P$ such that $f(x) = y$.
	\end{definition}
	
	\begin{theorem}
		Let $P$ be a poset and $\textswab{l}$ and $\textswab{r}$ be poset antiendomorphisms such that the following condition holds: \[\textswab{l}(\textswab{r}(\textswab{l}(x))) = \textswab{l}(x) \text{~and~} \textswab{r}(\textswab{l}(\textswab{r}(x))) = \textswab{r}(x).\] Then, $P$ satisfies the ascending chain condition for $\textswab{r}$-elements of $P$ if and only if $P$ satisfies the descending chain condition for $\textswab{l}$-elements of $P$.
	\end{theorem}
	
	\begin{proof}
		For the direct implication, let \[\textswab{l}(x_1) \geq \textswab{l}(x_2) \geq \cdots \geq \textswab{l}(x_n) \geq \cdots\] be a descending chain of $\textswab{l}$-elements in $P$. Since $\textswab{r}$ is order-reversing, we have the following ascending chain \[\textswab{r}(\textswab{l}(x_1)) \leq \textswab{r}(\textswab{l}(x_2)) \leq \cdots \leq \textswab{r}(\textswab{l}(x_n)) \leq \cdots\] which must stop somewhere by assumption. Therefore, there is a positive integer $m$ such that for any positive integer $i$ we have \[\textswab{r}(\textswab{l}(x_m)) = \textswab{r}(\textswab{l}(x_{m+i})).\] This implies that \[\textswab{l}(\textswab{r}(\textswab{l}(x_m))) = \textswab{l}(\textswab{r}(\textswab{l}(x_{m+i}))).\] Thus by properties of $\textswab{l}$ and $\textswab{r}$, we have \[\textswab{l}(x_m) = \textswab{l}(x_{m+i})\] showing $P$ satisfies the descending chain condition for $\textswab{l}$-elements of $P$. The converse implication is proved similarly and the proof is complete.
	\end{proof}
	
	Let $(P, \leq)$ be a poset and $A$ be a nonempty subset of $P$. It is said that $A$ is linearly ordered by $\leq$ if for all $x,y \in A$ either $x \leq y$ or $y \leq x$.
	
	\begin{theorem}
		Let $(P,\leq)$ be a poset and $\textswab{l}$ and $\textswab{r}$ be poset antiendomorphisms such that the following condition holds: \[\textswab{l}(\textswab{r}(\textswab{l}(x))) = \textswab{l}(x) \text{~and~} \textswab{r}(\textswab{l}(\textswab{r}(x))) = \textswab{r}(x).\] Then, $\textswab{r}$-elements of $P$ are linearly ordered by $\leq$ if and only if $\textswab{l}$-elements of $P$ are linearly ordered by $\leq$.
	\end{theorem}
	
	\begin{proof}
		Let $\textswab{r}$-elements of $P$ be linearly ordered by $\leq$. Imagine $\textswab{l}(x)$ and $\textswab{l}(y)$ are $\textswab{l}$-elements of $P$. Then, $\textswab{r}(\textswab{l}(x))$ and $\textswab{r}(\textswab{l}(y))$ are $\textswab{r}$-elements of $P$. By assumption \[\textswab{r}(\textswab{l}(x)) \leq \textswab{r}(\textswab{l}(y)) \text{~or~} \textswab{r}(\textswab{l}(y)) \leq \textswab{r}(\textswab{l}(x)).\] This implies that \[ \textswab{l}(\textswab{r}(\textswab{l}(x))) \geq \textswab{l}(\textswab{r}(\textswab{l}(y))) \text{~or~} \textswab{l}(\textswab{r}(\textswab{l}(y))) \geq \textswab{l}(\textswab{r}(\textswab{l}(x)))\] which implies that \[\textswab{l}(x) \geq \textswab{l}(y) \text{~or~} \textswab{l}(y) \geq \textswab{l}(x).\] Therefore, $\textswab{l}$-elements of $P$ are linearly ordered by $\leq$. The converse of the implication is proved similarly and the proof is complete.
	\end{proof}
	
	Let $X$ be a nonempty subset of a hemiring $H$. The left (right) annihilator of $X$ is defined as follows: \[\textswab{l}(X)= \{h\in H: hX = \{0\}\}~(\textswab{r}(X)= \{h\in H: Xh = \{0\}\}).\] Then, it is easy to see that the following properties hold:
	
	\begin{itemize}
		\item The left (right) annihilator of $X$ is a left (right) ideal of $S$.
		
		\item $X \subseteq \textswab{l}(\textswab{r} (X)) \cap \textswab{r}(\textswab{l} (X))$.
		
		\item If $X \subseteq Y \subseteq H$, then $\textswab{l}(Y) \subseteq \textswab{l}(X)$ and $\textswab{r}(Y) \subseteq \textswab{r}(X)$.
		
		\item $\textswab{l}(\textswab{r}(\textswab{l}(X))) = \textswab{l}(X)$ and $\textswab{r}(\textswab{l}(\textswab{r}(X))) = \textswab{r}(X)$.
	\end{itemize}
	
	In particular, we have the following:
	
	\begin{corollary}\label{orderlrannihilators}
		In any hemiring the a.c.c. for right annihilators is equivalent to the d.c.c. for left annihilators. On the other hand, right annihilators are linearly ordered by inclusion if and only if left annihilators are so.
	\end{corollary}
	
	Let $S$ be a subhemiring of a hemiring $H$. Note that if $I$ is a left annihilator of $S$, then there is a nonempty subset $X$ of $S$ such that \[I = \textswab{l}(X) = \{h\in S: hX =\{0\}\}.\] It is easy to see that $I$ is also an ideal of $H$ and so, a left annihilator of $H$. A similar statement holds for right annihilators. From this straightforward consideration, we have the following:
	
	\begin{proposition}
		Let $H$ be a hemiring. If $H$ satisfies a.c.c. or d.c.c. for left (right) annihilators of $H$, then so does any subhemiring of $H$.
	\end{proposition}
	
\section{Comparability of ideals of matrix semirings}\label{sec:chainconditionsmatrixsemirings}

Let us recall that for any semiring $S$ and positive integer $n$, if $I$ is an ideal of $S$ and then $M_n(I)$ is an ideal of the matrix semiring $M_n(S)$. On the other hand, any ideal of $M_n(S)$ is the form $M_n(I)$ for a unique ideal of $S$ \cite[Proposition 5.13]{Golan2003}.

\begin{lemma}\label{primenessmatrixsemiring}
If $P$ is a prime ideal of a semiring $S$, then $M_n(P)$ is a prime ideal of $M_n(S)$. On the other hand, if $\mathcal{P}$ is a prime ideal of $M_n(S)$, then there is a unique prime ideal $P$ of $S$ such that $\mathcal{P} = M_n(P)$. 
\end{lemma}

\begin{proof}
Let $I$ and $J$ be arbitrary ideals of $S$. It is easy to see that \[M_n(I)M_n(J) = M_n(IJ).\] Now, let $P$ be a prime ideal of $S$. Consider the ideal $M_n(P)$ of $M_n(S)$ and assume that $M_n(I) M_n(J) \subseteq M_n(P)$. This implies that $IJ \subseteq P$. Since $P$ is prime, we deduce that either $M_n(I) \subseteq M_n(P)$ or $M_n(J) \subseteq M_n(P)$. This shows that $M_n(P)$ is a prime ideal of $M_n(S)$.
	
On the other hand, let $\mathcal{P} = M_n(P)$ be a prime ideal of $M_n(S)$. Assume that $IJ \subseteq P$. Therefore, \[M_n(I)M_n(J) = M_n(IJ) \subseteq M_n(P).\] Primeness of $M_n(P)$ implies that either $I\subseteq P$ or $J\subseteq P$ and the proof is complete.	
\end{proof}

A ring is uniserial if its (two-sided) ideals are linearly ordered by inclusion \cite{Puninski2001}. Similarly, we define uniserial semirings as follows:
	
\begin{definition}
We say a semiring $S$ is uniserial if (two-sided) ideals of $S$ are linearly ordered by inclusion.
\end{definition}

Let us recall that a semiring $S$ is, by definition, proper if $S$ is not a ring \cite[p. 9]{HebischWeinert1998}.

\begin{example}
\label{uni-serial-semiring}
By Example 2.5 in \cite{Nasehpour2018}, the fuzzy semiring $(\mathbb I=[0,1], \max, \min)$ is a proper uniserial semiring. 
\end{example}
	
\begin{theorem}
Let $S$ be a semiring. Then, the following statements hold:

\begin{enumerate}
\item A semiring $S$ is uniserial if and only if the matrix semiring $M_n(S)$ is uniserial.
\item Prime ideals of a semiring $S$ are linearly ordered by inclusion if and only if prime ideals of $M_n(S)$ are so.
\item Subtractive (prime) ideals of $S$ are linearly ordered by inclusion if and only if subtractive ideals (prime) of $M_n(S)$ are so.
\end{enumerate}
\end{theorem}
	
\begin{proof}
(1): By Proposition 5.13 in \cite{Golan2003}, for each ideal $I$ of $S$, $M_n(I)$ is an ideal of $M_n(S)$ and each ideal of $M_n(S)$ is of the form $M_n(I)$ for some ideal $I$ of $S$. On the other hand, it is clear that $I \subseteq J$ if and only if $M_n(I) \subseteq M_n(J)$. Hence, $S$ is uniserial if and only if  $M_n(S)$ is so, as required. 	

(2): In view of Lemma \ref{primenessmatrixsemiring}, the proof is straightforward.

(3): Since an ideal $I$ of $S$ is subtractive if and only if $M_n(I)$ is subtractive, in view of the proof of (1) and Lemma \ref{primenessmatrixsemiring}, (3) holds and the proof is complete.
\end{proof}

\section{Pairs of commutative semirings with the same prime ideals}\label{sec:pairsofsemiringswiththesameprimeideals}

In this section, all semirings are commutative.

\begin{definition}\label{pseudovaluationsubsemiringdef} 
Let $S \subset U$ be two semidomains. We say $S$ is a pseudo-valuation subsemiring of $U$ if for each $u \in U \setminus S$ and non-unit $s\in S$, we have $u^{-1} s \in S$.
\end{definition}

\begin{lemma}\label{nonunitsthesame}
Let $S \subseteq T$ be two semirings. Then, the following statements hold:
\begin{enumerate}
	\item If $\Max(S) \subseteq \Spec(T)$, then a non-unit of $S$ is a non-unit of $T$.
	\item If $\Max(T) \subseteq \Spec(S)$, then a non-unit of $T$ is a non-unit of $S$. 
\end{enumerate}
\end{lemma}

\begin{proof}
(1): Let $s$ be a non-unit of $S$. So, there is a maximal ideal $\mathfrak{m}$ of $S$ such that $s \in \mathfrak{m}$. From the assumption $\Max(S) \subseteq \Spec(T)$, we obtain that $s$ is an element of a prime ideal of $T$. Since in any semiring, any prime ideal is a subset of a maximal ideal, we see that $s$ is an element of a maximal ideal of $T$ and so, it cannot be a unit of $T$. 

(2): Let $t$ be a non-unit in $T$. Then, $t$ is an element of a maximal ideal of $T$ and so by the assumption $\Max(T) \subseteq \Spec(S)$, $t$ is in a prime ideal of $S$. In particular, $t$ is an element of $S$ and cannot be a unit in $S$. This completes the proof.
\end{proof}

\begin{theorem}\label{pseudo-valuationsubsemiringthm} 
	Let $S \subseteq T \subset U$ be semidomains with $\Max(T) \subseteq \Spec(S)$ and $\Max(S) \subseteq \Spec(T)$. Then, $S$ is a pseudo-valuation subsemiring of $U$ if and only if $T$ is so.
\end{theorem}

\begin{proof}
	For the direct implication, assume that $u \in U \setminus T$ and $t \in T$ is a non-unit. By Lemma \ref{nonunitsthesame}, $t \in S$ cannot be a unit. Now, since by assumption $S$ is a pseudo-valuation subsemiring of $U$, we observe that $u^{-1}t \in S$. Finally, since $S \subseteq T$, we have $ u^{-1}t \in T$ showing that $T$ is also a pseudo-valuation subsemiring of $U$. Conversely, let $u \in U \setminus S$ and $s\in S$ be a non-unit. Consequently, by Lemma \ref{nonunitsthesame}, $s$ is not a unit in $T$. We distinguish two cases:
	
	Case I: If $u$ is not an element of $T$, then $u^{-1} s \in T$ because $T$ is a pseudo-valuation subsemiring of $U$ by assumption. Our claim is that $u^{-1} s$ is not a unit in $T$ since otherwise there is a $t\in T$ such that $(u^{-1}s)t = 1$ which implies that $u = st \in T$, a contradiction. In particular, by Lemma \ref{nonunitsthesame}, we have that $u^{-1}s \in S$.
	
	Case II: If $u$ is an element of $T$, then $u$ is a unit in $T$ since otherwise by Lemma \ref{nonunitsthesame}, $u$ will be an element of $S$, a contradiction. Our claim is that $u^{-1}s$ is a non-unit in $T$ since otherwise there is a $t\in T$ such that $u^{-1}s t = 1$ and so, $u = st$. Since $u \in T$ is a unit any factor of $u$ must be also a unit, a contradiction since by assumption $s$ is not a unit in $T$ (again by Lemma \ref{nonunitsthesame}). In particular, from $\Max(T) \subseteq \Spec(S)$, we obtain that $u^{-1}s \in S$.
	
	Hence, in any case, we see that $S$ is a pseudo-valuation subsemiring of $U$, as required.
	\end{proof}

The following is a generalization of Theorem 2.2 in \cite{KumarGaur2020}.

\begin{corollary}\label{pseudo-valuationsubsemiringcor}
Let $S \subseteq T \subset U$ be semidomains with $\Spec(S) = \Spec(T)$. Then, $S$ is a pseudo-valuation subsemiring of $U$ if and only if $T$ is so.
\end{corollary}

\section{Commutative semirings whose prime ideals are linearly ordered by inclusion}\label{sec:lopi}
	
	In this section, all semirings are commutative. Let us recall that if $\mathfrak{a}$ is an ideal of a semiring $S$, $\mathbb V(\mathfrak{a})$ is defined to be the set of all prime ideals $\mathfrak{p}$ of $S$ such that $\mathfrak{p} \supseteq \mathfrak{a}$. The radical of an ideal $\mathfrak{a}$, denoted by $\sqrt{\mathfrak{a}}$, is the intersection of all prime ideals containing it. So, $\sqrt{\mathfrak{a}} = \cap \mathbb V(\mathfrak{a})$. Krull's Theorem for semirings (see Proposition 7.28 in \cite{Golan1999(b)}) states that \[\sqrt{\mathfrak{a}} = \{s\in S: \exists ~n\in \mathbb N ~(s^n\in \mathfrak{a})\}. \]
	
	Similar to ring theory, we say an ideal $\mathfrak{a}$ of a semiring $S$ is radical if $\sqrt{\mathfrak{a}} = \mathfrak{a}$. For example, it is easy to see that if $\mathfrak{p}_1, \dots, \mathfrak{p}_n$ are prime ideals of $S$, then the ideal $\mathfrak{a} = \cap_{i=1}^{n} \mathfrak{p}_i$ is a radical ideal of $S$. The following is the semiring version of Theorem 1 in \cite{Badawi1995}:	
	
	\begin{theorem}
		
		\label{LOPmain}
		
		Let $S$ be a semiring. The following statements are equivalent:
		\begin{enumerate}
			\item The prime ideals of $S$ are linearly ordered.
			
			\item The radical ideals of $S$ are linearly ordered.
			
			\item Each proper radical ideal of $S$ is prime.
			
			\item The radicals of the principal ideals of $S$ are linearly ordered.
			
			\item For each $x,y \in S$, there is an $n \geq 1$ such that either $x|y^n$ or $y|x^n$.
		\end{enumerate}

		\begin{proof}
			$(1) \Rightarrow (2)$: Suppose that $\mathfrak{a}$ is a proper ideal of $S$. By Krull's Theorem for semirings, the radical of an ideal equals the intersection of all prime ideals containing it. Now, since the prime ideals of $S$ are linearly ordered, the intersection of all the primes containing $\mathfrak{a}$ is the minimum prime among all primes of $S$ over $\mathfrak{a}$.
			
			$(2) \Rightarrow (3)$: Let $\mathfrak{a}$ be a proper ideal of $S$. Similar to the proof of the implication $(1) \Rightarrow (2)$, we see that the radical $\sqrt{\mathfrak{a}}$ of the ideal $\mathfrak{a}$ equals to the minimum prime ideal over $\mathfrak{a}$.
			
			$(3) \Rightarrow (1)$: Let $\mathfrak{p},\mathfrak{q}$ be two distinct prime ideals of $S$ and consider the ideal $\mathfrak{a}=\mathfrak{p} \cap \mathfrak{q}$. Therefore, by Proposition 7.30 in \cite{Golan1999(b)}, $\sqrt{\mathfrak{a}}=\mathfrak{a}$, and so by assumption, $\mathfrak{a}$ is a prime ideal of $S$. Now, since $ \mathfrak{p} \mathfrak{q} \subseteq \mathfrak{p} \cap \mathfrak{q}$, it is obvious that either $\mathfrak{p} \subset \mathfrak{q}$ or $\mathfrak{q} \subset \mathfrak{p}$. 
			
			The statement $(2) \Rightarrow (4)$ requires no comment.
			
			$(4) \Rightarrow (5)$: Let $x,y\in S$. By assumption, either $\sqrt{(x)} \subseteq \sqrt{(y)}$ or $\sqrt{(y)} \subseteq \sqrt{(x)}$. Therefore, either $x^k \in (y)$ for some $k\in \mathbb N$ or $y^l \in (x)$ for some $l\in \mathbb N$. Now, let $n=\max\{k,l\}$. It is clear that either $x|y^n$ or $y|x^n$.
			
			$(5) \Rightarrow (1)$: Let $\mathfrak{p},\mathfrak{q}$ be two distinct prime ideals of $S$ with $x \in \mathfrak{p}\setminus \mathfrak{q}$. Therefore, for every $y \in \mathfrak{q}$ there is an $n \geq 1$ such that $x|y^n$. This implies that $y\in \mathfrak{p}$.
		\end{proof}
		
	\end{theorem}
	
	The proof of the following is straightforward:
	
	\begin{proposition}
	Let $S$ be a local semiring of Krull dimension 1. Then, prime ideals of $S$ are linearly ordered by inclusion.
	\end{proposition}
	
	Let us recall that $a$ and $b$ in a semiring $S$ are associates if $a=ub$ for some unit $u\in U(S)$. A semidomain $S$ is said to be a GCD semidomain if $\gcd(a,b)$ exists for any $a,b\in S$, whenever at least one of the elements $a$ and $b$ is nonzero \cite[Definition 4.4]{Nasehpour2020}. Also, note that a semidomain $S$ is a valuation semiring if and only if for any $a$ and $b$ in $S$ either $a$ divides $b$ or $b$ divides $a$ \cite[Theorem 2.4]{Nasehpour2018}. W.V. Vasconcelos shows that if $A$ is a commutative domain, then $A$ is a valuation domain if and only if it is a GCD domain and the prime ideals are linearly ordered \cite[Proposition A]{Vasconcelos1972}. Now we prove its semiring version as follows: 
	
	\begin{theorem}
		
		\label{GCDvaluationvasconcelos}
		
		Let $S$ be a semidomain. Then, $S$ is a valuation semiring if and only if $S$ is a GCD semidomain which has prime ideals that are linearly ordered.
		
		\begin{proof}
			
			It is clear, by the note above, that if $S$ is a valuation semiring, then $S$ is a GCD semidomain. Now
			suppose that $x,y$ are two nonzero nonunit elements of $S$ and $z=\gcd(x,y)$ and suppose that $z$ is associated to neither $x$ nor $y$. Let $a=x/z$ and $b=y/z$. Then, neither $a$ nor $b$ is a unit of $S$. Therefore, by Theorem \ref{LOPmain}, there exists an $m \geq 1$ such that either $a|b^m$ or $b|a^m$. Clearly, $\gcd(a,b)=1$. Also by \cite[Proposition 4.6(3)]{Nasehpour2020}, we can conclude that for every $n \geq 1$, $\gcd(a,b^n)=\gcd(b,a^n)=1$. Therefore, $a$ or $b$ is a unit of $S$, a contradiction. Hence, either $x|y$ or $y|x$. This finishes the proof.
		\end{proof}
	\end{theorem}
	
	Let us recall that a prime ideal $\mathfrak{p}$ of an integral domain $D$ is called strongly prime, in the sense of Hedstrom and Houston \cite{HedstromHouston1978}, if whenever $x,y \in K$ and $xy \in \mathfrak{p}$, then $x \in \mathfrak{p}$ or $y \in \mathfrak{p}$, where $K$ is the field of fractions of the integral domain $D$. If every prime ideal of $D$ is strongly prime, then $D$ is called a pseudo-valuation domain [abbreviated as PVD]. Inspired by these concepts, we give the following definitions:
	
	\begin{definition}
		
		\label{PVS}
		
		Let $S$ be a semidomain.
		
		\begin{enumerate}
			\item We define a prime ideal $\mathfrak{p}$ of $S$ to be strongly prime if whenever $x,y \in F(S)$ and $xy \in \mathfrak{p}$, then $x \in \mathfrak{p}$ or $y \in \mathfrak{p}$.
			
			\item We define a semidomain $S$ to be a pseudo-valuation semidomain (abbreviated as PVS) if every prime ideal of $S$ is strongly prime.
		\end{enumerate}
		
	\end{definition}
	
	\begin{proposition}
		
		\label{valuationPVS}
		
		Every valuation semiring is a PVS.

		\begin{proof}
			Let $S$ be a valuation semiring and let $\mathfrak{p}$ be a prime ideal of $S$. Suppose that $xy \in \mathfrak{p}$ where $x,y \in F(S)$, the semifield of fraction of $S$. If both $x$ and $y$ are in $S$, we are done. Suppose that $x \notin S$. Since $S$ is a valuation semiring, we have $x^{-1} \in S$ \cite[Theorem 2.4]{Nasehpour2018}. Hence, $y=xyx^{-1} \in \mathfrak{p}$. This finishes the proof.
		\end{proof}
	\end{proposition}
	
	\begin{remark}

		Let $V$ be a valuation domain of the form \[V=K+M :=\{k+m | k\in K, m\in M\},\] where $K \subseteq V$ is a field and $M$ is the unique maximal ideal of $V$, e.g. $V=K[[X]]$. Now, let $F$ be a proper subfield of $K$. Then, it is easy to verify that \[R=F+M :=\{f+m | f\in F, m\in M\}\] is a subring of $V$ which is a pseudo-valuation domain but not a valuation one (see Example 2.1 in \cite{HedstromHouston1978}). The following question is still open for us:
		
		\begin{question}
			
			Is there any proper semiring which is a PVS but not a valuation one?
		\end{question}
		
	\end{remark}
	
	\begin{proposition}
		
		\label{PVSLOP}
		
		Let $S$ be a PVS. Then, the prime ideals of $S$ are linearly ordered.
		\begin{proof}
			
			Let $S$ be a PVS. Also, let $\mathfrak{p}$ be a prime and $\mathfrak{a}$ an arbitrary ideal of $S$. Our claim is that $\mathfrak{p}$ and $\mathfrak{a}$ are comparable. If not, then there exist an $x\in \mathfrak{p}$ and $s\in \mathfrak{a}$ such that $x \notin \mathfrak{a}$ and $s\notin \mathfrak{p}$. Note that $(x/s)\cdot s =x \in \mathfrak{p}$. Also, since $x/s \notin S$, $x/s \notin \mathfrak{p}$. On the other hand, $s \notin  \mathfrak{p}$, contradicting this assumption that $\mathfrak{p}$ is strongly prime. Therefore, all prime ideals of $S$ are linearly ordered and the proof is complete.
		\end{proof}
	\end{proposition}
	
	\begin{theorem}
		
		\label{PVSchar1} 
		
		A semidomain $S$ is a PVS if and only if $S$ is quasi-local and the only maximal ideal $\mathfrak{m}$ of the semiring $S$ is strongly prime.

		\begin{proof}
			
			It is clear that $\Rightarrow $ holds. Now, we prove the inverse as follows:
			
			$(\Leftarrow)$: Suppose that the only maximal ideal $\mathfrak{m}$ of the quasi-local semidomain $S$ is strongly prime. Let $\mathfrak{a}$ be a prime ideal of $S$ and $x,y \in F(S)$ and $xy \in \mathfrak{a}$. If $x,y \in S$, then $x \in \mathfrak{a}$ or $y \in \mathfrak{a}$. Now, suppose $x \notin S$. Since $xy \in \mathfrak{m}$ and $x \notin S$, we have $y \in \mathfrak{m}$. Suppose $y \notin \mathfrak{a}$. Then, $y^2 \notin \mathfrak{a}$ and so, $d=(y^2/xy) \notin S$. Note that $dx=y \in \mathfrak{m}$ while $x$ and $d$ are not elements of $\mathfrak{m}$, a contradiction. Thus, $y \in \mathfrak{a}$ and $\mathfrak{a}$ is strongly prime.
		\end{proof}
		
	\end{theorem}
	
	\begin{lemma}
		
		\label{xinverse}
		
		Let $S$ be a PVS. Then, for every $x \in F(S)\setminus S$, $x^{-1} \mathfrak{p} \subseteq \mathfrak{p}$, where $\mathfrak{p}$ is a prime ideal of $S$.

		\begin{proof}
			Since $S$ is a PVS, $\mathfrak{p}$ is strongly prime. If $x \in F(S)\setminus S$ and $y \in \mathfrak{p}$, then $yx^{-1}x=y \in \mathfrak{p}$. Thus either $yx^{-1} \in \mathfrak{p}$ or $x \in \mathfrak{p}$. Since $x \notin S$, we must have $yx^{-1} \in \mathfrak{p}$. Hence, $x^{-1} \mathfrak{p} \subseteq \mathfrak{p}$.
		\end{proof}
	\end{lemma}
	
	Let us recall that a semiring $S$ is quasi-local if and only if $S \setminus U(S)$ is an ideal of $S$, where $U(S)$ is the set of all unit elements of $S$ (see Example 6.1 and Proposition 6.61 in \cite{Golan1999(b)}). The following result is a semiring version of Proposition 3.1 in \cite{Anderson1979}:
	
	\begin{theorem}
		
		\label{subsetlocal}
		
		Let $(S,\mathfrak{m})$ be a quasi-local semidomain. Then, $S$ is a PVS if and only if for every $x \in F(S)$, either $xS \subseteq \mathfrak{m}$ or $\mathfrak{m} \subseteq xS$.

		\begin{proof}
			$(\Rightarrow)$: If $x \in S$, then either $xS \subseteq \mathfrak{m}$ or $\mathfrak{m} \subseteq xS$, depending on whether $x$ is a nonunit or unit of $S$. If $x \in F(S) \setminus S$, then $x^{-1} \mathfrak{m} \subseteq \mathfrak{m}$ by Lemma \ref{xinverse}. Thus, $\mathfrak{m} \subseteq x \mathfrak{m} \subseteq xS$.
			
			$(\Leftarrow)$: We show that $\mathfrak{m}$ is strongly prime. Let $xy \in \mathfrak{m}$ with $x,y \in F(S)$. Assume $x \notin \mathfrak{m}$, so $\mathfrak{m} \subseteq xS$. Therefore, $xy \in xS$, so $y \in S$. If $y \notin \mathfrak{m}$, then $y$ is a unit of $S$ and so, $y^{-1} \in S$. This implies $x \in \mathfrak{m}$, a contradiction. Hence, $y \in \mathfrak{m}$. \end{proof} 
	\end{theorem}	
	
	\begin{theorem}
		
		\label{PVSchar2}
		
		Let $N=S\setminus U(S)$ be the set of all nonunit elements of a semidomain $S$. Then, the following statements are equivalent:
		
		\begin{enumerate}
			
			\item The semiring $S$ is a PVS with the maximal ideal $N$.
			\item For each pair of $\mathfrak{a}, \mathfrak{b}$ of ideals of $S$, either $\mathfrak{b} \subseteq \mathfrak{a}$ or $\mathfrak{a} \cdot \mathfrak{c} \subseteq \mathfrak{b}$ for every proper ideal $\mathfrak{c}$ of $S$.
			\item For every $x,y \in S$, either $yS \subseteq xS$ or $xzS \subseteq yS$ for every nonunit $z \in S$.
			\item For every $x,y \in S$, either $x|y$ or $y|xz$ for every nonunit $z \in S$.
			\item For every $x,y \in S$, either $yS \subseteq xS$ or $xN \subseteq yS$.
			\item For every $x,y \in S$, either $yN \subseteq xS$ or $xS \subseteq yN$. 
		\end{enumerate}
		
		\begin{proof}
			
			(1) $\Rightarrow$ (2): Let $\mathfrak{a}, \mathfrak{b}$ be ideals of $S$ and $\mathfrak{c}$ be a proper ideal of $S$. Suppose that $\mathfrak{b}$ is not a subset of $\mathfrak{a}$ and $\mathfrak{a} \cdot \mathfrak{c}$ is not a subset of $\mathfrak{b}$. So, there exists $x \in \mathfrak{b\setminus a}$ and $yz \in \mathfrak{a} \cdot \mathfrak{c}$ for some $y \in \mathfrak{a}$ and $z \in \mathfrak{c}$ such that ${x/y} \in F(S) \setminus S$ and $yz/x \in F(S) \setminus S$. On the other hand, $(x/y)(zy/x)=z \in \mathfrak{c} \subseteq N$ and neither $x/y \in N$ nor $yz/x \in N$, which is a contradiction. 
			The implications (2) $\Rightarrow$ (3) and (3) $\Rightarrow$ (4) are obvious. 
			
			(4) $\Rightarrow$ (1): Suppose that for every $x,y \in S$ and any nonunit $z \in S$, either $x|y$ or $y|xz$. Let $x$ be a nonunit element of $S$ and $y \in S$. Then, either $x|y$ or $y|x^2$. Hence, the prime ideals of $S$ are linearly ordered by Theorem \ref{LOPmain}. In particular, $S$ is quasi-local with the maximal ideal $N$ \cite[Proposition 3.5]{Nasehpour2018}. By Theorem \ref{PVSchar1}, we just need to show that $N$ is strongly prime. Suppose that $uv \in N$ for some $u,v \in F(S)$. If $u \in S$ or $v \in S$, then it is easy to see that $u \in N$ or $v \in N$. Therefore, suppose that $u,v \in F(S) \setminus S$. Let $u=y/x$ and $v=z/w$ for some $x,y,z,w \in S$. Since $u=y/x \in F(S) \setminus S$ and $uv=yz/xw \in N$, $y|x(yz/xw)$. Thus, $v=z/w \in S$, a contradiction. Hence, if $uv \in N$ for some $u,v \in F(S)$, then $u \in N$ or $v \in N$. 
			
			It is clear that the statements (4) and (5) are equivalent. Now, we prove that (6) implies (1).
			
			$(6) \Rightarrow (1)$: Let $x,y \in S$. It is obvious that (6) implies (5). On the other hand, (5) and (4) are equivalent. Therefore, similar to the proof of the implication (4) $\Rightarrow$ (1), we have that either $x|y$ or $y|x^2$. Therefore, the prime ideals of $S$ are linearly ordered. In particular, $S$ is quasi-local with the maximal ideal $N$. Hence, $S$ is a PVS  by Theorem \ref{subsetlocal}. Finally, the statement $(1) \Rightarrow (6)$ is just a restatement of Theorem \ref{subsetlocal}.  
		\end{proof}
		
	\end{theorem}
	
	\begin{theorem}
		
		\label{LOP}
		
		Let $(S,\mathfrak{m})$ be a quasi-local semidomain. If for each $x \in F(S)$, $x \mathfrak{m} \subseteq \mathfrak{m}$ or $\mathfrak{m} \subseteq x \mathfrak{m}$, then the set of prime ideals of $S$ is linearly ordered.	
		
		\begin{proof}
			Let $\mathfrak{p}$ and $\mathfrak{q}$ be two distinct prime ideals of $S$ with $\mathfrak{p} \not\subseteq \mathfrak{q}$. Choose $x \in \mathfrak{p\setminus q}$. Our claim is that $xy^{-1} \mathfrak{m}$ cannot be a subset of $\mathfrak{m}$ for any $y \in \mathfrak{q}$. In the contrary, take $y \in \mathfrak{q}$ such that $xy^{-1} \mathfrak{m} \subseteq \mathfrak{m}$. So, $x \mathfrak{m} \subseteq y \mathfrak{m} \subseteq \mathfrak{q}$. Since $\mathfrak{q}$ is prime, either $\mathfrak{m} \subseteq \mathfrak{q}$ or $x \in \mathfrak{q}$. The latter condition contradicts choice for $x$, as $x \in \mathfrak{p\setminus q}$. Now, if $\mathfrak{m} \subseteq \mathfrak{q}$, by maximality of $\mathfrak{m}$, $\mathfrak{m} = \mathfrak{q}$ and since $x\in \mathfrak{p} \subseteq \mathfrak{m}$, we have $x \in \mathfrak{m} = \mathfrak{q}$, again a contradiction. Thus $yx^{-1} \mathfrak{m} \subseteq \mathfrak{m}$, so $y \mathfrak{m} \subseteq \mathfrak{p}$. Therefore, either $y \in \mathfrak{p}$ or $\mathfrak{m} \subseteq \mathfrak{p}$. Now, if $\mathfrak{m} \subseteq \mathfrak{p}$, then $\mathfrak{m} = \mathfrak{p}$, which implies that $y\in \mathfrak{p}$. Therefore, $\mathfrak{q} \subseteq \mathfrak{p}$ and the proof is complete.
		\end{proof}
	\end{theorem}
	
	Let us recall that a prime ideal $\mathfrak{p}$ of a semiring $S$ is a divided prime ideal of $S$ if $\mathfrak{p} \subset (x)$ for every $x \in S \setminus \mathfrak{p}$. A semiring $S$ is divided if each prime ideal of $S$ is divided \cite[Definition 1.3]{BehzadipourNasehpour2020}.
	
	\begin{theorem}
		
		\label{dividedchar}
		
		The following statements for a semidomain $S$ are equivalent:
		
		\begin{enumerate}
			
			\item The semiring $S$ is divided.
			
			\item For every pair of ideals $\mathfrak{a},\mathfrak{b}$ of $S$, the ideals $\mathfrak{a}$ and $\sqrt{\mathfrak{b}}$ are comparable.
			
			\item For every $x,y \in S$, the ideals $(x)$ and $\sqrt{(y)}$ are comparable.
			
			\item For every $x,y \in S$, either $x|y$ or $y|x^n$  for some $n \geq 1$. 
			
		\end{enumerate}

		\begin{proof}
			
			$(1) \Rightarrow (2)$: Let $S$ be a divided semiring and $\mathfrak{a},\mathfrak{b}$ two proper ideals of $S$. Since $S$ is divided, the prime ideals of $S$ are linearly ordered \cite[Proposition 1.5]{BehzadipourNasehpour2020}. Also by Theorem \ref{LOPmain}, we can conclude that $\sqrt{\mathfrak{b}}=\mathfrak{p}$ is prime. So by \cite[Proposition 1.5]{BehzadipourNasehpour2020}, the ideals $\mathfrak{a}, \sqrt{\mathfrak{b}}$ are comparable. 
			
			The statements $(2) \Rightarrow (3)$ and $(3) \Rightarrow (4)$ are clear. 
			
			$(4) \Rightarrow (1)$: Suppose that for every $x,y \in S$, there is a natural number $n \geq 1$ such that either $x|y^n$ or $y|x$. Let $\mathfrak{p}$ be a prime ideal of $S$ and $w \in S \setminus \mathfrak{p}$ and $z \in \mathfrak{p}$. Since $z$ does not divide $w^n$ for any $n \geq 1$, it follows that $w|z$. Therefore, $\mathfrak{p}$ is comparable to every principal ideal of $S$. Hence, $S$ is a divided semidomain.
		\end{proof}
	\end{theorem}
	
	\begin{corollary}
		
		\label{dividedLO}
		
		Let $S$ be a divided semidomain. Then, the prime ideals of $S$ are linearly ordered by inclusion.
		
		\begin{proof} It follows from Theorem \ref{LOPmain}.\end{proof}
	\end{corollary}
	
	\subsection*{Acknowledgments}
	
	During the research and composition of the current paper, the second author has been supported by the Department of Basic Sciences at the Golpayegan College of Engineering and his special thanks go to the Department for providing all the necessary facilities available to him for successfully conducting this research.

\end{document}